\documentclass{article}
\usepackage{amsmath,amssymb,amsthm,latexsym}

\numberwithin{equation}{section}

\theoremstyle{definition}
\newtheorem{theorem}{Theorem}

\newtheorem{definition}[theorem]{Definition}
\usepackage{authblk}
\newcommand{\V}{\mathsf{V}}
\newcommand{\En}{E_{n+1}}
\newcommand{\EE}{E_{2\cdot 3^{2n+1}+2}}
\newcommand{\1}{1^\text{'}}

\begin{document}
\title{On the equational complexity of RRA}
\author{Jeremy F. Alm}
\affil{Department of Mathematics\\ Illinois College\\ 1101 W. College Ave.\\Jacksonville, IL 62650\\ \texttt{alm.academic@gmail.com}}
\date{\today}
% prints the title
\maketitle
% suppresses page numbers
\thispagestyle{empty}
% body of text

\begin{abstract}
We prove that the equational complexity function for the variety of
representable relation algebras is bounded below by a log-log
function.
\end{abstract}

\section{Introduction}

Let \textsf{RRA} denote the class of representable relation
algebras. \textsf{RRA} is definable by equations \cite{Tar55}, but
not by finitely many \cite{Mon64}.  Indeed, any equational basis
must contain equations containing arbitrarily many variables
\cite{Jon91}.  It is an open question whether \textsf{RRA} is
definable by first-order formulas using some bounded number of
variables---see \cite{HH}, page 625.

A \emph{weak representation} of a relation algebra  is an
isomorphism to an \textsf{RRA} that doesn't necessarily preserve the
operations $+$ and $-$ but does preserve $\cap$.  Let \textsf{wRRA}
denote the class of weakly representable relation algebras.
\textsf{wRRA} is not finitely based \cite{HM00}, and \textsf{RRA} is
not finitely based over \textsf{wRRA} \cite{And94}. It was recently
shown that \textsf{wRRA} is a variety \cite{Pes09}.  Since
\textsf{RRA} has no finite-variable equational basis it must be the
case that at least one of the following holds:

\begin{enumerate}
\item  \textsf{wRRA} has no finite-variable equational basis;
\item there is no finite-variable equational basis that defines \textsf{RRA} over
\textsf{wRRA}.
\end{enumerate}

It would be interesting to know which of these hold.  The author
submits this to the reader as an open problem.

All of these results speak to the ``bad behavior" of  \textsf{RRA} .
In this note, we want to focus on a related question for finite
algebras: given a finite $A\in \textsf{RA}$, how much of the
equational theory of \textsf{RRA} do we have to verify  in $A$
before we know that $A\in \textsf{RRA}$?

\section{Definitions}

We take the following definition from \cite{MSW08}:

\begin{definition}
The \emph{length} of an equation is the total number of operation
symbols and variables appearing in the equation.  For a variety $\V$
of finite signature, the \emph{equational complexity} of $\V$ is
defined to be a function $\beta_\V$ such that for a positive integer
$m$, $\beta_\V(m)$ is the least integer $N$ such that for any
algebra $A$ of the similarity class of $\V$ with $|A|\leq m$, $A\in
\V$ iff $A$ satisfies all equations true in $\V$ of length at most
$N$.
\end{definition}

 For example, the
length of $(x + y)\cdot z = x\cdot z + y\cdot z$  is 12. We note
that for a variety $\V$ of finite signature, $\beta_\V$ always
exists. To see this, fix $m$, and consider the collection of
algebras in the similarity class of $\V$ of size of most $m$ that
are not in $\V$. For each algebra in the collection, take the
shortest equation that witnesses the algebra's non-membership in
$\V$.  Let $\ell$ be the length of the longest such shortest
equation. Then $\ell +1$ is an upper bound for $\beta_\V(m)$.

Throughout the rest of this paper, let $\V= \textsf{RRA} $. In
\cite{Lyn61}, Roger Lyndon gave a general construction of relation
algebras from projective geometries.  We are interested in the
algebras that come from finite projective lines, and we will use
them to find a lower bound on $\beta_\V$. We give a definition here
that is equivalent to the one Lyndon gave.

Let $\En$ be a finite integral relation algebra with $n$ symmetric
diversity atoms $a_1, \ldots, a_n$  and one identity atom $\1$.
Composition on the atoms is defined thus:  $$a_i ; a_i = \1 + a_i
\quad\text{ and }\quad a_i ; a_j = \overline{a_i + a_j + \1} \text{
for } i\neq j$$

Lyndon proved that $\En$ is representable iff there exists a
projective plane of order $n-1$.  Bruck and Ryser proved in
\cite{BruRys49} that there is no projective plane of order $2 \cdot
3^{2n+1}$; hence, $\EE$ is non-representable.  However, every proper
subalgebra $A$ of $\En$ embeds into $E_{p+1}$ for any prime $p>n$,
and hence is representable.  J\'{o}nsson used this fact in
\cite{Jon91} to give a proof that $\mathsf{RRA}$ has no $k$-variable
basis for $k<\omega$.  This implies that $\beta_\V(m)$ is not
bounded above.

\section{The lower bound}

The computation of this lower bound follows the proof of Lemma  6 in
\cite{MSW08}. Consider $E_{2\cdot 3^{2n+1}+2}$: since there is no
projective plane of order $2\cdot 3^{2n+1}$, $E_{2\cdot
3^{2n+1}+2}\not\in\textsf{RRA}$. Therefore, there is some equation
$\varepsilon$ such that $\textsf{RRA}\models\varepsilon$ but
$E_{2\cdot 3^{2n+1}+2}\not\models\varepsilon$.  We recall that every
proper subalgebra of $E_{2\cdot 3^{2n+1}+2}$ is representable.
Consider the number of distinct variables in $\varepsilon$, and
suppose that it is no more than $k\leq\log_2 3 \cdot (2n+1)$. Then
take $b_1, \ldots, b_k \in E_{2\cdot 3^{2n+1}+2}$.  The subalgebra
generated by $b_1, \ldots, b_k$ is the boolean subalgebra generated
by $\1, b_1, \ldots, b_k$.  This subalgebra is no larger than
$2^{2^{k+1}}$, and thus is proper, since
$$2^{2^{\log_2 3 \cdot (2n+1) +1}} = 2^{2\cdot 3^{2n+1}} < 2^{2\cdot
3^{2n+1} + 2} = |E_{2\cdot 3^{2n+1}+2}|$$

Thus we can conclude that $\varepsilon$ contains more than $\log_2 3
\cdot (2n+1)$ variables, since any equation of fewer variables true
in all representable relation algebras would have to be satisfied by
the representable subalgebra of $E_{2\cdot 3^{2n+1}+2}$ generated by
$b_1, \ldots, b_k$. Now consider the length of $\varepsilon$: since
$\varepsilon$ contains $k$ distinct variables, it must contain at
least $k-2$ binary operation symbols, hence its length is at least
$2k-2$. This gives us that $$2\log_2 3 \cdot (2n+1) -2 < \beta_\V
\left( 2^{2\cdot 3^{2n+1} + 2} \right) \qquad (\star)$$

Now choose $m\in \mathbb{Z}^+$, with $m\geq 2^8$. Then there is some
$n\in\mathbb{Z}^+$ so that $$2^{2\cdot 3^{2n+1} + 2 } \leq m \leq
2^{2\cdot 3^{2n+3} + 2 }$$  Then $m \leq 2^{2\cdot 3^{2n+3} + 2 }$ \
gives us that   $$\frac{1}{2}\log_3\left(\frac{1}{2}\log_2 (m) -
1\right) - \frac{3}{2} \leq n \qquad (\star\star)$$  Let $f(n) =
2\log_2 3 \cdot (2n+1) -2$.  We apply $f$ to both sides of
$(\star\star)$, which (since $f$ is increasing) yields
\begin{align*}
  2\log_2 3 \cdot (\log_3 \left(\frac{1}{2} \log_2 (m) -1\right) -2) -2 &\leq 2\log_2 3 \cdot (2n+1) -2 \\
                                                                       &< \beta_\V \left( 2^{2\cdot 3^{2n+1}+2}\right) && \text{by $(\star)$}\\
                                                                       &\leq \beta_\V (m),
\end{align*}
\noindent where the last line follows from the monotonicity of
$\beta_\V$.\\

Therefore $\beta_\V(m) >2\log_2 3 \cdot (\log_3 \left(\frac{1}{2} \log_2 (m) -1\right) -2) -2$ for all $m\geq 2^8$.  \\

Since the size of a finite relation algebra is always a power of 2,
we can make some aesthetic changes. Let $M$ be the number of atoms
of a finite algebra $A$, and let $\beta^*_\V$ be the equational
complexity function that takes as input the number of atoms of an
algebra (rather than the cardinality).  Then we get $$\beta^*_\V (M)
> 2\log_2 3 \cdot[\log_3 (M/2-1) -2] -2$$

\section{Conclusion}

Since the language of \textsf{RA} has finite signature, $\beta_\V$
is always finite. In \cite{MSW08}, tools are given for finding upper
bounds for locally finite varieties.  For the variety \textsf{RRA},
the derivation of an upper bound may prove more difficult.  The
author submits this as another open problem.

\end{document}